\documentclass[12pt, a4paper]{amsart}
\input xy
\xyoption{all}
\numberwithin{equation}{section}

\theoremstyle{plain}
  \begingroup
        \newtheorem{theorem}[equation]{Theorem}
        
        \newtheorem{proposition}[equation]{Proposition}
        \newtheorem{corollary}[equation]{Corollary}
        
        \newtheorem{assumption}[equation]{Assumption}
        \newtheorem{remark}[equation]{Remark}
  \endgroup

\theoremstyle{definition}
  \begingroup
        \newtheorem{definition}[equation]{Definition}

        \newtheorem{notation}[equation]{Notation}
  \endgroup

\newcommand{\diagrama}{\xymatrix}
\newcommand{\mr}[1]{\buildrel {#1} \over \longrightarrow}

\begin{document}


\title{LOCALIC GALOIS THEORY} 

\author{Eduardo J. Dubuc}

\maketitle


\begin{abstract}
In this article we prove the following:

\emph{A topos with a point is 
connected atomic if and 
only if it is the classifying topos of a localic group, 
and this group can be taken to be the locale of automorphisms of the 
point}. 

We explain and give the necessary definitions to understand this 
statement.

The hard direction in this equivalence was first proved in print  
in \cite{JT}, Theorem 1, Section 3, Chapter VIII, and it follows from 
a characterization of atomic topoi in terms of open maps and from
a theory of descent for morphisms of topoi and locales.

We develop our version and our proof of this theorem, which is completely 
independent of descent theory and of any other result in \cite{JT}.
Here the theorem follows as an straightforward consequence of a direct 
generalization of the fundamental theorem of Galois.

In Proposition I of 
``Memoire sur les conditions de resolubilite
des equations par radicaux'', Galois established that any 
intermediate extension of the splitting field of a polynomial with 
rational coefficients is the fixed field of its galois group. 
We first state and prove the (dual) categorical interpretation of of this 
statement, which is a theorem about atomic sites with a \emph{representable} 
point. These developments  correspond exactly to 
\emph{Classical Galois Theory}. 

In the general case, the point determines a 
proobject  and it becomes (tautologically) prorepresentable. 
 We state and prove the, mutatus mutatis, prorepresentable 
version of Galois theorem. In this 
case the \emph{classical} group of automorphisms has to be 
replaced by the \emph{localic} group of automorphisms.
These developments form the content of 
a theory that we call \emph{Localic Galois Theory}.

\end{abstract}

\vspace{4ex}

 INTRODUCTION 

\vspace{2ex}

In this article we prove the following:

\emph{Theorem B: A topos $\mathcal{E}$
with a point $\mathcal{E}ns \mr{p} \mathcal{E}$, $p^{*} = F$, is 
connected atomic if and 
only if it is the classifying topos $\mathcal{B}G$ of a localic group $G$, 
and this group can be taken to be $lAut(F) = lAut(p)^{op}$}. 

We shall 
explain and give the necessary definitions to understand this 
statement.

The hard direction in this equivalence was first proved in print  
in \cite{JT}, Theorem 1, Section 3, Chapter VIII, and it follows from 
the characterization of atomic topoi as those topoi such that  
 $\mathcal{E} \rightarrow \mathcal{E}ns$ and the diagonal
 $\mathcal{E} \mr{\Delta} \mathcal{E} \times \mathcal{E}$ are open, 
 and from a theory of descent for morphisms of topoi and locales.

We develop our version and our proof of this theorem, which is completely 
independent of descent theory and of any other result in \cite{JT}.
Here the theorem follows as an straightforward consequence of a direct 
generalization of the fundamental theorem of Galois.

In Proposition I of 
``Memoire sur les conditions de resolubilite
des equations par radicaux'', (see \cite{E}), Galois established that any 
intermediate extension of the splitting field A of a polynomial with 
rational coefficients is the fixed field of its galois group. 
We first state and prove the (dual) categorical interpretation (Theorem \ref{Galois} 
below) of of this statement (in the dual category the fixed field is 
the quotient of A by the action of the galois group). This 
interpretation is a theorem about atomic sites. From this result 
(\ref{Galois}) it follows 
in an straightforward manner the following:

\emph{Theorem A: A topos $\mathcal{E}$
with a representable point $\mathcal{E}ns \mr{p} \mathcal{E}$, 
\mbox{$p^{*} = [A,\;-],\; A \in \mathcal{E}$} is 
connected atomic if and 
only if it is the classifying topos $\mathcal{B}G$ of a discrete group $G$, 
and this group can be taken to be $Aut(A)^{op}$}.

We call this development the
\emph{representable case} of Galois Theory, and it corresponds exactly to 
\emph{Classical Galois Theory}. 

In the general case, the point $p^{*}$ determines a 
proobject $P$ and it becomes (tautologically) prorepresentable 
$p^{*} = [P,\;-]$.
Theorem  \ref{Galois2} below
is just the prorepresentable version of Theorem \ref{Galois}. In this 
case the classical group of automorphisms $Aut(P)$ has to be 
replaced by the localic group $lAut(P)$ (whose points form the group
$Aut(P)$ which can be trivial).
Theorem B follows from this result (\ref{Galois2}) in the same straightforward 
manner than Theorem A from \ref{Galois}.

We call this development the
\emph{prorepresentable case} of Galois Theory, and it is the content of 
a theory that we call \emph{Localic Galois Theory}.

\vspace{2ex}

A \emph{localic space} is the formal dual of a local, and 
a \emph{localic group} is a group object in the category of localic 
spaces.

 Our basic approach 
is to work with locales, and consider them as 
posets (as in topos theory one works with topoi as categories). After all, 
one can hardly expect to prove all the results about localic 
spaces which are false for topological spaces if one works 
``as if they were topological spaces'' (arguments justified by test maps 
and the like). It is 
some times surprising how localic techniques are very often more 
simple and clear than its dual geometrical counterparts. Geometrical 
intuition and experience is important to an overall understanding, 
but when it comes to prove basic results it is of little help (however 
it is useful, as in topos theory, in order to develop implication chains 
utilizing basic results).

\vspace{2ex}

This paper is divided in eight short sections, with the principal 
contributions in sections 2., 4. and 6.

\vspace{2ex}

 1. Classical Galois Theory.
 
 2. The theorems of localic Galois Theory.
 
 3. Preliminaries on localic spaces and groups.
 
 4. The locale of automorphisms of a set valued functor.
 
 5. The (pre) topology generated by a family of covers.
 
 6. Proof of the theorems of localic Galois theory.
 
 7. Preliminaries on the classifying topos of a localic group.
 
 8. Characterization of the classifying topos of a localic group. 

\vspace{2ex}
 
In section 1 we state and prove our interpretation of the theorems 
of classical Galois Theory, (theorems \ref{transitivity} and 
\ref{Galois}), and then we prove theorem A in the introduction.

\vspace{2ex}

In section 2 we state explicitly the theorems of localic Galois Theory
(theorems \ref{transitivity2} and \ref{Galois2}).

\vspace{2ex}
 
In section 3 we fix the terminology and notation on locales, 
and recall some necessary facts.

\vspace{2ex}

In section 4 we develop a fundamental construction in this paper, 
namely, that of the localic group $lAut(F)$ of automorphisms of a set valued 
functor $\mathcal{C} \mr{F} \mathcal{E}ns$. To prove theorems 
\ref{transitivity2} and  \ref{Galois2} the 
straightforward construction, as the appropriate subspace of the product 
$\prod_{X \in \mathcal{C}} \: lAut(FX)$, is useless.
 
Unlike the locale $lAut(X)$ of automorphisms of a set $X$, the locale of
 relations $lRel(X)$ 
is  functorial on $X$, with values in the 
category of Posets. This is technically of great importance, since it 
allows to develop constructions corresponding to Grothendieck's 
construction of the (co)-fibered category associated to a functor with 
values in the category of Categories. We exploit this in our 
construction of the locale $lAut(F)$, for a set valued functor 
$F$, \emph{by first constructing the locale $lRel(F)$}, and then 
the locale (subspace) of 2-valued sheaves for the 
 Grothendieck topology that forces a relation to be a bijection. 
 We were inspired by Gavin Wraith presentation in \cite{W} of 
 the locale of automorphisms of a set X.

\vspace{2ex}
 
Section 5 is technical on the generation of grothendieck topologies 
out of some basic covers.

\vspace{2ex}

Section 6 (together with section 4) contains the important 
contributions made in this paper. Here we prove 
the fundamental theorems of localic Galois theory. That is, 
the theorems of section 2 on the prorepresentable case (theorems 
\ref{transitivity2} and \ref{Galois2}).

\vspace{2ex}     

In section 7 we recall several necessary
facts on the category of sets furnished with a continuous action of a 
localic group. Although this facts 
are widely believed to be true, nobody has cared to prove them in print.

\vspace{2ex}

Finally, in section 8 we show how theorem B follow from the theorems 
in section 2.
       
\vspace{4ex} 

\section{Classical Galois Theory}    

This corresponds to the representable case of the theory. Notice that 
in the category dual of the category of intermediate extensions of the 
splitting field $A$ of a polynomial with rational coefficients, the fixed 
field of a group $H \subset Aut(A)$ is, categorically, the quotient of A by 
the action of $H$.
 
Let $\mathcal{C}$ be any category and $A \in \mathcal{C}$ be any object. 
Assume:

\begin{assumption} \label{assumption} 
\ \newline \indent
   i) Every arrow $Y \mr{} X$ in $\mathcal{C}$ is an strict 
   epimorphism.

  ii) For every $X \in \mathcal{C}$ there exists $A \mr{} X$.

 iii) The representable functor $F = [A, -]$ preserves strict epimorphisms.

\end{assumption}

Then:

\begin{theorem} \label{transitivity}
For every object $X \in \mathcal{C}$ the action of the group 
$Aut(A)^{op}$ on the set $[A, X]$ is transitive.  
\end{theorem}

\begin{theorem}[Galois Theorem] \label{Galois}
Every arrow $A \mr{x} X$  in $\mathcal{C}$ is the 
categorical quotient of $A$ by the action of the Galois group
 \mbox{$Fix(x) = \{ h \in Aut(A) \: | \: xh = x \} \, \subset \, Aut(A)$.}
\end{theorem}

These theorems follow easily from the following proposition:

\begin{proposition} \label{basica} 
Every arrow 
$X \mr{f} A$ is an isomorphism. In particular, every endomorphism of $A$ is 
an isomorphism, $Aut(A) = [A, A]$.
\end{proposition}

\begin{proof}
In fact, 
from iii) it follows that there is $A \mr{g} X$
such that $fg = id$. Then, $g$ is a monomorphism. Since by i) it is also an 
strict epimorphism, it follows that it is an isomorphism, and consequently so 
is $f$.
\end{proof}


{\bfseries  Proof of theorems \ref{transitivity} and \ref{Galois}}  
\begin{proof}
Theorem \ref{transitivity} follows immediately from iii). Let now $A \mr{x} X$,
 and  assume $A \mr{y} Y$ is any arrow such 
that $Fix(x) \subset Fix(y)$. Since $x$ is an strict epimorphism
 (see \ref{strictepi}), to prove theorem \ref{Galois} it will be 
enough to show that given any two arrows
$\diagrama{ Z \ar@<1ex>[r]^{t} \ar@<-1ex>[r]^{s} & A }$ 
, the implication \mbox{``$\, xs = xt \; \Longrightarrow \; ys = yt$\,''} holds. 
By \ref{basica} we can assume $Z = A$ \ and $s$ 
invertible. Let $xs = xt$, then $ts^{-1} \in Fix(x)$, thus also 
$ts^{-1} \in Fix(y)$. Thus $ys = yt$.
\end{proof}

\vspace{2ex}

An straightforward consequence of \ref{transitivity} and 
\ref{Galois} is theorem A:

\begin{theorem}
 A topos $\mathcal{E}$
with a representable point $\mathcal{E}ns \mr{p} \mathcal{E}$, 
\mbox{$p^{*} = [A,\;-],\; A \in \mathcal{E}$} is 
connected atomic if and 
only if it is the classifying topos $\mathcal{B}G$ of a discrete group $G$, 
and this group can be taken to be $Aut(A)^{op}$.
\end{theorem}
\begin{proof}
 By theorem \ref{transitivity} the functor $[A, -]$ lifts into the category 
 of transitive $G$-sets, for $G \;=\; [A, A]^{op} \;=\; Aut(A)^{op}$ 
(\ref{basica}). Theorem \ref{Galois} then essentially means that this 
lifting is full and faithful. Since every transitive $G$-set is a 
quotient of the $G$-set $[A, X] \;for\; X = A$, it follows 
(by the comparison lemma (\cite{G2}, 
Expose III, 4.) that the topos of sheaves for the canonical 
topology on $\mathcal{C}$ is equivalent to the topos of $G$-sets. 

It is immediate to check (see \ref{atomicsite})
that the data in assumption \ref{assumption} is a 
connected atomic site with a                              
representable point, and any connected atomic topos with a representable 
point can be presented in this way (see \cite{BD}). This finishes the proof.
\end{proof}

\vspace{2ex}

\section{The theorems of localic Galois Theory}
\vspace{2ex}

This corresponds to the prorepresentable case of the theory.

Let $\mathcal{C}$ be any category and $P \in Pro \,\mathcal{C}$ be any 
pro-object (in the sense of Grothendieck \cite{G2}). 
Recall that the proobject $P$ is the formal dual of a prorepresentable 
functor $F:\mathcal{C} \mr{} \mathcal{E}ns$, $F = [P, -]$. We shall 
call the functor $F$ to be the \emph{fiber functor}. 
 
Assume:

\begin{assumption} \label{assumption2}
\ \newline \indent
   i) Every arrow $Y \mr{} X$ in $\mathcal{C}$ is an strict 
   epimorphism.

  ii) For every $X \in \mathcal{C}$ there exists $P \mr{} X$. That 
  is, $FX \neq \emptyset$.

 iii) The prorepresentable functor $F = [P, -]$ preserves strict epimorphisms.

\end{assumption}

 Then:

\begin{theorem} \label{transitivity2}
For every object $X \in \mathcal{C}$ the action of the localic group 
of automorphisms $lAut(P)^{op} = lAut(F)$ on the set $[P, X] = FX$ is transitive.  
\end{theorem}

\begin{theorem} \label{Galois2}
Every arrow $P \mr{x} X$ in $\mathcal{C}$ is the 
categorical quotient, relative to the category $\mathcal{C}$, of $P$ by the 
action of the Galois group $lFix(x) \in lAut(P)$ described 
informally as $\{h \in Aut(P) \: | \: xh = x \}$.
 
The equivalent version of this statement in terms of the fiber functor 
is reminiscent of the lifting lemma in classical covering theory:

Lifting Lemma: Given any objects 
$X \in \mathcal{C},\; Y \in \mathcal{C}$, and elements  
$x \in FX,\; y \in FY$, 
if $lFix(x) \:\leq\: lFix(y)$ in $lAut(F)$, then there exist a unique 
arrow \mbox{$X \mr{f} Y$} in $\mathcal{C}$ such that $F(f)(x) = y$.
  
\end{theorem}

In the rest of this paper we shall explain and prove these theorems, 
giving the necessary definitions.

\vspace{2ex}

\section{Preliminaries on Localic Spaces and Groups} 

Topoi are often considered as generalized topological spaces, but the intuition 
in topos theory is not only geometrical. We think of locale theory as a 
reflection of topos theory (with the poset $2 = \{0, 1\}$ playing the 
role of the category $\mathcal{E}ns$ of sets), as well as that of a theory of 
generalized topological spaces.

\vspace{2ex}
 
We consider a \emph{poset} as a category, and in this vein a 
\emph{partial order} is a reflexive and transitive 
relation, not necesarially antisymetric. We denote the order relation 
either by $'' \rightarrow''$ or by $''\leq''$.  We shall call \emph{objects} 
the elements of a poset.

A morphism of posets is an 
\emph{injection} when it is injective in the isomorphism classes. That is, 
if it creates isomorphisms.


\vspace{2ex}

A \emph{locale} is a complete lattice in which finite infima 
distribute over arbitrary suprema. A morphism of locales 
$E \mr{f^{*}} H$ 
is defined to be a function $f^{*}$ preserving finite infima and 
arbitrary suprema (notice that we put automatically an upper star 
to indicate that these arrows are to be considered as 
inverse images of geometric maps). We shall also refer to such a 
morphism as \emph{an $H$-valued point of $E$}. 2-valued points 
$E \mr{f^{*}} 2$ are just called \emph{points}.

\vspace{2ex}

 \emph{Inf-lattices} $D$ are sites of definition for locales (rather than bases 
 of opens). 2-valued \emph{presheaves} $D^{op} \rightarrow 2$ correspond to
  downward closed
  subsets $T$, and they form a locale, $D^{\wedge} = 2^{D^{op}}$. Given a 
  Grothendieck 
 (pre) topology on $D$, 2-valued \emph{sheaves} are those T such that 
 for each cover $u_{\alpha} \rightarrow u \,,\; (\forall \alpha \; 
 u_{\alpha} \in T) \;\; \Rightarrow \;\; (u \in T)$, and they also form a 
 locale, denoted $D^{\sim}$. The associated sheaf defines a morphism of locales 
  \mbox{$D^{\wedge} \rightarrow  D^{\sim}$,} and this is a procedure in which 
  quotients of locales are obtained. A \emph{site} is in this sense a 
  \emph{presentation} of the locale of sheaves. 

  An \emph{$H$-valued point} of an inf-lattice is an 
 inf-preserving morphism into a locale $H$. When $H = 2$, it 
corresponds to an upward closed subset $P$ such that 
\mbox{$u \in P \,,\: w \in P \;\Rightarrow\; u \wedge w \in P$}.
 An $H$-valued point of a site in addition must send covers into 
epimorphic families. When $H = 2$, this corresponds to the usual requirement. 
 That is, for each cover \mbox{$u_{\alpha} \rightarrow u$, $(u \in P) \;
  \Rightarrow \; (\exists \, \alpha \; u_{\alpha} \in P)$}. 

 The basic result of this construction is that the associate sheaf
\mbox{$D \mr{\#} D^{\sim}$} is a point which is generic, in the sense 
that giving any locale $H$, composing with $\#$ defines an equivalence of posets
\mbox{$\mathcal{P}oints(D^{\sim}, H) \mr{\simeq} \mathcal{P}oints(D, H)$}.
 \emph{Points of a site of definition and of the locale of sheaves are 
the same thing}.

\vspace{2ex}

 A \emph{localic space} is the formal 
dual of a local. Thus, $E \mr{f^{*}} H$ defines a map or morphism of 
localic spaces from $H$ to $E$, $H \mr{f} E$. Following \cite{JT}, all 
these maps are called \emph{continuous maps}.  
 A \emph{point} of a localic space $E$ is a point of the corresponding 
 locale.

 A function preserving 
finite infima is an injection if and only if it reflects 
isomorphisms (recall that finite infima determines the order relation).    

A \emph{surjection} between localic 
spaces is a map whose inverse image reflects isomorphisms.
 A locale has \emph{enough points} if its family of points is 
(collectively) surjective.

Any local $E$ determines a topology (in the classical sense) on its 
set of points by means of the correspondence, for $u \in E$, 
$P \in \mathcal{P}oints(E)$ and $U \subset \mathcal{P}oints(E)\,$:   
$\;\;P \in U\;\Leftrightarrow\; u \in P$.     

A localic space is a (sober) topological space if and only if it has 
enough points. In this case, the topology in $\mathcal{P}oints(E)$ also 
determines $E$ since $u \simeq v \;\Leftrightarrow\; U = V$.

\vspace{2ex}

A \emph{localic monoid}, (resp. \emph{localic group}) is a monoid 
object (resp. group object) in the category of localic spaces. 
A \emph{morphism of monoids (or groups)} $H \mr{\varphi} G$ is a continuous 
map such that 
\mbox{$m^{*}\varphi^{*} = (\varphi^{*}\otimes \varphi^{*})\,m^{*}$}
(where $m$ denotes the multiplication in the two structures).

\vspace{2ex}

We recall now a construction of the free inf-lattice on 
a poset $D$. That is, the inf-completion of $D$, which is the inf-lattice 
$\mathcal{D}(D)$ whose points $\mathcal{D}(D) \rightarrow H$ correspond 
exactly to the order preserving morphisms $D \rightarrow H$. Warning: 
these are not the points of $D$ !. 
 
 
\begin{proposition}  \label{freeinflattice}
Given any poset $D$, consider the diagram:
$$
\diagrama@1
        {
         D \;\; \ar@<+2pt> `u[r] `[rr] ^-{Yoneda} [rr] 
                        \ar @{^{(}->}[r]   
          & \;\; \mathcal{D}(D) \;\;  \ar @{^{(}->}[r]
          & \;\; (2^{D})^{op}
        }
$$ 
where $\mathcal{D}(D) \subset (2^{D})^{op}$ is the full subposet of finitely 
generated upward closed subsets of $D$. Given a finite subset 
\mbox{$\{a_{1}, \ldots ,\,a_{n}\} \subset D$}, We denote 
\mbox{$[A] \;=\; [<a_{1}>, \ldots ,\,<a_{n}>] \;=\; \{a \in D \;|\; 
\exists i \; a_{i} \leq a\}$}. If \mbox{$[B] = \,[<b_{1}>, \ldots ,\,<b_{k}>]$}, 
then $[A] \rightarrow [B]$ in $\mathcal{D}(D)$ (that is 
 \mbox{$[B] \subset [A]$)} if and only if 
\mbox{$\;\exists\:\sigma\,:\,  \{1, \ldots, k\} 
\mr{} \{1, \ldots, n\},\; a_{\sigma i} \leq b_{i}$}\,.       
\end{proposition}

A particular case of this construction is the free inf-lattice on a 
set $X$,
$\mathcal{D}(X)$, which is the poset of 
finite subsets of $X$ with the reverse of the natural order.
It follows that the 
free locale on $X$ is the locale of presheaves on $\mathcal{D}(X)$, 
$\mathcal{L}(X) = \mathcal{D}(X)^{\wedge}$.

The points of $\mathcal{L}(X)$  are (by definition) the subsets of $X$. 
If $x \in X$, we denote
 $[<\!x\!>]$ the corresponding generator in $\mathcal{L}(X)$. 
 If $S \subset X$, we have $S \in [<\!x\!>]\; \Leftrightarrow\; x \in S$.
Similarly, if $\{x_{1}, \ldots ,\,x_{n}\} \subset X$, we 
write \mbox{$[<x_{1}>, \ldots ,\,<x_{n}>]$} for the 
corresponding object in $\mathcal{D}(X) \subset \mathcal{L}(X)$. 
Notice that this 
object defines the open set (in the topological space of points)   
\mbox{$\{S \subset X \:|\: x_{i} \in S,\: i = 1, \ldots ,n\}$}. The 
following is clear:  

\begin{proposition}  \label{lRel(X)}
The locale of relations $lRel(X)$ on a set X is the free locale 
$\mathcal{L}(X \times X) = \mathcal{D}(X \times X)^{\wedge}$. If 
$\{(x_{1},\,y_{1})\,\ldots , (x_{n},\,y_{n})\} \subset X \times X$, we write 
$[<x_{1}\,|\,y_{1}>\,\ldots, <x_{n}\,|\,y_{n}>]$ for the corresponding 
object in the site or in the locale. 
\end{proposition}     

\vspace{2ex}

We take now from \cite{W} a site of definition for the localic group 
of automorphisms of a set. That is, a localic group such that its 
points are the automorphisms of X.  

\begin{proposition}  \label{lAut(X)}
The locale of automorphisms $lAut(X)$ on a set $X$ is the locale of 
sheaves on the site with underline poset the inf-lattice 
$\mathcal{D}(X \times X)$, and with the covers generated by 
the following families (in the notation in \ref{lRel(X)}):: 
$$
\begin{array}{l} 
 \emptyset \;\rightarrow\; 
[<\!z\,|\,x\!>,\;<\!z\,|\,y\!>]\,,  \\
 \emptyset \;\rightarrow\; 
[<\!x\,|\,z\!>,\;<\!y\,|\,z\!>]\,,  \\
 \;\;\;\;\;\;\;\;\;\;\;\;\;\;\;\;\;\;\;\;\;\;\;\;\;\;\;\;\;\;\;\;\;\;\;\;\;\;\;
          (each \; x,\,y,\,z, \: x \neq y)  \\
\,[<\!x\,|\,z\!>] \;\rightarrow \; 1, \; x \in X\,,   \\
\,[<\!z\,|\,x\!>] \;\rightarrow \; 1, \; x \in X\,,   \\
\;\;\;\;\;\;\;\;\;\;\;\;\;\;\;\;\;\;\;\;\;\;\;\;\;\;\;\;\;\;\;\;\;\;\;\;\;\;\;
           (each \, z)\,.
\end{array}
$$            
\end{proposition}
\begin{proof}
It follows immediately from \ref{lRel(X)}. The coverings above force a 
relation to be, in turn, univalued, injective, everywhere defined, and surjective.
\end{proof}
 
We shall abuse the notation and omit to indicate the associate sheaf 
morphism $lRel(X) \rightarrow lAut(X)$. Thus, 
\mbox{$[<x_{1}\,|\,y_{1}>\,\ldots, <x_{n}\,|\,y_{n}>]$} also 
denotes the corresponding object in $lAut(x)$.

Actually, this locale has enough points, 
and it is the usual set of bijections of X furnished with the product 
topology. We have then the usual open set in the base of this topology
\mbox{$\{f:\, X \rightarrow X \;|\; f(x_{i}) = y_{i} \: i = 1, \ldots ,n\}$}. 

 The motivation in G. Wraith paper was to consider this 
presentation in an arbitrary topos, where it defines a local which 
in general will not have enough points.

The local $lRel(X)$ is a localic monoid, and its binary operation 
restricts to $lAut(X)$ and defines a localic group. This structure is 
given by:
$$
m^{*}([<\!x \:|\:y\!>]) \;=\;
 \bigvee\nolimits_{z} \;[<\!x \:|\:z\!>]\otimes[<\!z \:|\:y\!>]\;\;\; 
$$
The identity map $X \rightarrow X$ determines a point 
$lAut(X) \mr{e^{*}} 2$ which is the neutral element for m:
 $e^{*}[<\!x\,|\,y\!>] \;=\; 1 \;\;\Leftrightarrow\;\; x \;=\; y$. Thus 
$e \in [<\!x\,|\,y\!>] \;\;\Leftrightarrow\;\; x \;=\; y$.  

All this is described in \cite{W}, from were we take also the 
definition of \emph{action} of 
a localic group $G$ on a set $X$, see \ref{action} below.

\vspace{2ex}

\section{The locale of automorphisms of a set-valued functor}

 Recall that given any category $\mathcal{C}$ and any functor 
 $\mathcal{C} \mr{F} \mathcal{E}ns$, the diagram of $F$, 
that we denote $\Gamma_{F}$, is the category whose objects are
the elements of the disjoint union of the sets $FX, \; X \in 
\mathcal{C}$. That is, pairs $(x,X)$ where $x \in FX$. The arrows 
$(x,X) \mr{f} (y,Y)$ are maps $X \mr{f} Y$ such that $F(f)(x)=y$.
There is a diagram
$$
\diagrama@R=0pt
  {
   {\Gamma_{F}}^{^{op}}\ar[r] & \mathcal{E}ns^{\, \mathcal{C}} \\
   (x, X) \ar@{|->}[r]  &  [X, -]
                           }
$$
with the obvious definition on arrows, and $F$ is the colimit of
this diagram.



\begin{definition} \label{deF}
We define a poset, that we denote $D_{F}$, 
by the following rule:
$$
  \frac{(x,\:X) \leq (y,\:Y)}
  {\; \exists \; X \mr{f} Y \; F(f)(x) = y \;}
$$ 
\end{definition}

Consider the function $FX \mr{\lambda_{X}} D_{F}$ defined by 
$\lambda_{X}(x) = (X,\;x)$. 
The proof of the following proposition is immediate:

\begin{proposition}  \label{genericparts}
The poset  $D_{F}$ has, and therefore it is characterized, by the 
following universal property: 

For each $X \in \mathcal{C}$, there is a function 
$FX \mr{\lambda_{X}} D_{F}$, and for each $X \mr{f} Y$ a transformation 
$\lambda_{X} \rightarrow \lambda_{Y} \circ F(f)$ (that is, for each 
\mbox{$x \in FX,\;$} $\lambda_{X}(x) \: \leq \: (\lambda_{Y} \circ F(f))(x)$).
 And for any other 
such data, there is a unique morphisms of posets $\phi$ (as indicated in the  
diagram below) such that $\phi \circ \lambda_{X}\:=\: \phi_{X}\,,\;
\phi \circ \lambda_{Y}\:=\: \phi_{Y}\,$:
$$
\diagrama
        {
         FX \ar[rd]^{\lambda_{X}}  \ar@(r, ul)[rrd]^{\phi_{X}} \ar[dd]_{F(f)}  \\
         & \;\;D_{F}\;\; \ar@{-->}[r]^{\phi} & \;D  \\
         FY \ar[ru]^{\lambda_{Y}} \ar@(r, dl)[rru]^{\phi_{Y}} 
        } 
$$ 
\end{proposition}

\begin{definition} \label{natural relation}
A \emph{natural relation} is a 
relation $R \subset F \times F$ in the functor category. That is, it is a 
family of relations $RX$ on $FX$, $X \in \mathcal{C}$, such that given 
any arrow $X \mr{f} Y$ in $\mathcal{C}$, and
 $(x_{0}, x_{1}) \in FX \times FX$:
$$ (x_{0}, x_{1}) \in RX \;\; \Rightarrow \;\; (F(f)(x_{0}), 
\,F(f)(x_{1})) \, \in RY$$
In other terms, it is a family of functions $FX \times FX 
\mr{\phi_{X}} 2$ such that 
$$\phi_{X}(x_{0}, x_{1}) \:
 \leq \: (\phi_{Y} \circ (F(f) \times F(f))(x_{0}, x_{1})$$
It is clear that if a natural relation is functional, then it is a 
natural transformation.
\end{definition}  

\vspace{2ex}

 Consider the composite of the diagonal functor 
 $\mathcal{C} \rightarrow \mathcal{C} \times \mathcal{C}$
  with $F \times F$, that we denote $\Delta F$, 
$(\Delta F)(X) \,=\, FX \times FX$. Notice that there is a full 
and faithful inclusion of categories 
$\Gamma_{F} \hookrightarrow  \Gamma_{\Delta F}$ and consequently a full 
inclusion of posets $D_{F} \hookrightarrow  D_{\Delta F}$.

From proposition  
\ref{genericparts} it follows immediately (see also definition 
\ref{deF}):

\begin{proposition}[Generic Natural Relation]  \label{genericnr}
The poset  $D_{\Delta F}$ has, and therefore it is characterized, by the 
following universal property: 

For each $X \in \mathcal{C}$, there is a function 
$FX \times FX \mr{\lambda_{X}} D_{\Delta F}$, and for each $X \mr{f} Y$ a 
transformation 
$\lambda_{X} \rightarrow \lambda_{Y} \circ (F(f) \times F(f))$. That is, 
for each $(x_{0}, x_{1}) \in FX \times FX$, 
$$\lambda_{X}(x_{0}, x_{1}) 
\: \leq \: (\lambda_{Y} \circ (F(f) \times F(f))(x_{0}, x_{1})$$

 And for any other 
such data, there is a unique morphisms of posets $\phi$ (as indicated in the  
diagram below) such that $\phi \circ \lambda_{X}\:=\: \phi_{X}\,,\;
\phi \circ \lambda_{Y}\:=\: \phi_{Y}\,$:

 
$$
\diagrama
        {
         FX \times FX \ar[rd]^{\lambda_{X}}  \ar@(r, ul)[rrd]^{\phi_{X}} 
                                             \ar[dd]_{F(f) \times F(f)}  \\
         & \;\;D_{\Delta F}\;\; \ar@{-->}[r]^{\phi} & \;D  \\
         FY \times FY \ar[ru]^{\lambda_{Y}} \ar@(r, dl)[rru]^{\phi_{Y}} 
        } 
$$

It follows that a morphisms of posets $D_{\Delta F} \rightarrow 2$ 
corresponds exactly to the data defining a natural relation of $F$.  
\end{proposition}

\begin{corollary} \label{lRel(F)}
The points of the locale of presheaves 
$\mathcal{D}(D_{\Delta F})^{\wedge}$ on the free inf-lattice 
$\mathcal{D}(D_{\Delta F})$ on the poset  $D_{\Delta F}$ 
are exactly the natural relations of $F$. That is, 
$lRel(F) \,=\, \mathcal{D}(D_{\Delta F})^{\wedge}$    
\end{corollary}

For later reference, and according with \ref{lRel(X)} and 
\ref{freeinflattice}, we record:

\begin{notation} \label{freeinflattice2}
By definition, the set of objects of $D_{\Delta F}$ is the disjoint 
union of the sets $FX \times FX,\; X \in \mathcal{C}$. Given an 
element $(X,\: (x_{0},\,x_{1}))$ and a finite subset $A \subset D_{\Delta F}$
of this set, we denote  
$$[(X,\,<\!x_{0}\,|\,x_{1}\!>),\; A]
\;=\; [\{(X,\: (x_{0},\,x_{1}))\} \cup A] \; = \;
[(X,\,<\!x_{0}\,|\,x_{1}\!>)]\,\wedge\, [A]$$ 
the corresponding object in $\mathcal{D}(D_{\Delta F})$.
\end{notation}

We shall construct now the locale of automorphisms of a set valued functor 
$F$ by defining a site structure on the inf-lattice $\mathcal{D}(D_{\Delta F})$.

\begin{proposition}  \label{lAut(F)}
The locale of automorphisms $lAut(F)$ of a set-valued functor $F$
is the locale of sheaves on the site with underline poset the 
 inf-lattice  $\mathcal{D}(D_{\Delta F})$, and with the covers generated by 
the following families:
$$
\begin{array}{l} 
 \emptyset \;\rightarrow\; 
[(X,\;<\!z\,|\,x\!>),\;(X,\;<\!z\,|\,y\!>)]\,,  \\
 \emptyset \;\rightarrow\; 
[(X,\;<\!x\,|\,z\!>),\;(X,\;<\!y\,|\,z\!>)]\,,  \\
 \;\;\;\;\;\;\;\;\;\;\;\;\;\;\;\;\;\;\;\;\;\;\;\;\;\;\;\;\;\;\;\;\;\;\;\;\;\;\;
          (each\; X,\; and \; each \; x \neq y,\,z\, \in FX)  \\
\,[(X,\;<\!x\,|\,z\!>)] \;\rightarrow \; 1, \; x \in FX\,,   \\
\,[(X,\;<\!z\,|\,x\!>)] \;\rightarrow \; 1, \; x \in FX\,,   \\
\;\;\;\;\;\;\;\;\;\;\;\;\;\;\;\;\;\;\;\;\;\;\;\;\;\;\;\;\;\;\;\;\;\;\;\;\;\;\;
           (each \;X \;and \;each \; z \in FX)\,.
\end{array}
$$            

Recall that the object 
$[(X_{1},\, <\!x_{1}\,|\,y_{1}\!>)\,\ldots, (X_{n},\, 
  <\!x_{n}\,|\,y_{n}\!>)]$
corresponds to the open set  
$\{\phi:\, F \rightarrow F \;|\; \phi X_{i}(x_{i}) = y_{i}\}$
\end{proposition}
\begin{proof}
It follows immediately from \ref{lAut(X)} and \ref{lRel(F)}.
\end{proof}

For each $X \in \mathcal{C}$ the map
$FX \times FX \mr{\lambda_{X}} \mathcal{D}(D_{\Delta F})$, 
defined by $\lambda_{X}^{*}[<\!x_{0}\,|\,x_{1}\!>] = [(X,\,<\!x_{0}\,|\,x_{1}\!>)]$
determines a morphism of locales $lRel(F) \mr{\lambda_{X}} lRel(FX)$.
 
The locale $lRel(F)$ is a localic monoid with the binary operation 
$m$ defined by 
$m^{*}[(X, \;<\!x \:|\: y\:>] = \lambda_{X}^{*}m_{X}^{*}[<\!x \:|\: 
y\:>]$, and by definition the arrow $\lambda_{X}$ becomes a morphism of
 monoids.
This operation on $lRel(F)$ restricts and defines the group structure 
of $lAut(F)$.

Consider now the morphism of locales given by the associate sheaf
\mbox{$\mathcal{D}(D_{\Delta F}) \hookrightarrow lRel(F) \mr{\#} lAut(F)$}.
 We have:
\begin{proposition} \label{actionF}
For each object $X \in \mathcal{C}$, the composite of the maps: 
$$FX \times FX \mr{\lambda_{X}} \mathcal{D}(D_{\Delta F}) \mr{\#} lAut(F)$$
 $(\#\,\lambda_{X})([<\!x_{0} \,|\, x_{1}\!>]) \; = \; 
 \#[(X,\; <\!x_{0} \,|\, x_{1}\!>)]$,
determines a morphism of locales $lAut(FX) \rightarrow lAut(F)$ that 
defines an action of $lAut(F)$ on the set $FX$. Furthermore, given 
any arrow $X \mr{f} Y$, the function $FX \mr{F(f)} FY$ becomes a 
morphism of actions (see section 6).  
\end{proposition}


\begin{proof}
For the first assertion it suffices to show that this map sends covers into 
covers on the respective
 sites of definition. But this is clear. The second assertion follows 
 by the diagram in \ref{genericnr}.
\end{proof}

\begin{remark} \label{identidad}
The identity $F \rightarrow F$ determines a point $lAut(F) \mr{e^{*}} 2$ 
given by $e^{*}[(X,\;<\!x\,|\,z\!>)] = 1  \;\Leftrightarrow\; x = z$. 
This point is the neutral element for the group structure. We have
  \mbox{$e \in [(X,\;<\!x\,|\,z\!>)] \;\Leftrightarrow\; x = z$}.    
\end{remark}     

\vspace{2ex}

\section{The (pre)topology generated by a family of covers}

Suppose we have a category $\mathcal{D}$ (with finite limits to 
simplify), and a family of \emph{(basic) covers}
$D_{\alpha} \rightarrow D$ on some objects $D \in \mathcal{D}$. We 
consider the site determined by the (pre)-topology generated by these 
covers. To check that a finite limit preserving functor is a point for 
this site it is enough to test the point condition only on the basic 
covers. To check that a presheaf is a sheaf, it is enough to test 
the sheaf condition on all covers obtained by pulling-back basic 
covers (in some terminology, the \emph{covering system} generated by 
the basic covers). However, in this paper we have to deal with a more 
subtle problem. We have to check that a given presheaf $T$ (which is not a 
sheaf) behaves as a sheaf against some given objects $A \in \mathcal{D}$. 
In this case, it is necessary to test the sheaf condition on all the covers 
$A_{\alpha} \rightarrow A$ of the (pre)-topology generated. This is 
so because $T$ may not be a sheaf against the objects $A_{\alpha}$, 
and this fact breaks the argument used to show that it is enough to test 
the sheaf condition on the covering system.

We need a careful description (by transfinite induction) of the (pre) 
topology generated by a family of basic covers.

First we shall fix some notation. Let $Cov$ be a collection of 
\emph{small} families 
$A_{\alpha} \rightarrow A \in Cov(A)$ 
 of arrows (to be considered as 
\emph{coverings}) on each object $A \in \mathcal{D}$ (small in the sense that the 
index $\alpha$ ranges over a set in $\mathcal{E}ns$). 
Given a collection $Cov$, define a new collection, denoted $\pi Cov$:

$$
   A_{\alpha} \rightarrow A \; \in \; \pi Cov(A) \;\; 
   \Longleftrightarrow \;\; \exists \;
    B_{\alpha} \rightarrow B \; \in \; Cov(B)
$$
$$
   \;\;\;\;\;\;\;\;\;\;\;\;\;\;\;\;\;\;\;\;\;\;\;\;\;\;\;\;\;\;\;\;\;\;
   \;\;\;\;\;\;
   and \; a \; pullback \;\; 
   \diagrama
     {
       A_{\alpha} \ar[r] \ar[d]  &  A \ar[d]  \\
       B_{\alpha} \ar[r]  &  B 
     }
$$
Given another collection $Dov$, define the \emph{composite} 
$Dov \, \circ \, Cov$ by means of the following implication:
$$
  \begin{array}{l}
  A_{\alpha} \rightarrow A \; \in \; Cov(A) \;\; and \;\;
  \forall \, \alpha \;\; A_{\alpha,\, \beta} \rightarrow A_{\alpha} \; \in \; 
  Dov(A_{\alpha})  \\
  \;\;\;\;\;\;\;\;\;\;\;\;\;\;\;\;\;\;\;\;\;\;\;\;\;\;\;\;\;
  \;\;\;\;
  \Longrightarrow \;\;  
  A_{\alpha, \, \beta} \rightarrow A_{\alpha} \rightarrow A \; 
  \in \; (Dov \circ Cov)(A) 
  \end{array}
$$
Notice that collections compose the other way than arrows, and that 
the two constructions above preserve the size condition. Let $Iso$ 
be the collection whose covers consists of a single isomorphism.

\vspace{2ex}

A \emph{covering system} is a collection $Cov$ such that 
$Iso \subset Cov$ and  
$\pi Cov \subset Cov$. A covering system is a \emph{(pre) topology} 
if in addition $Cov \circ Cov \subset Cov$.


\begin{proposition} \label{pregenerated}
Let $\mathcal{D}$ be a category with finite limits, and 
\mbox{$D_{\alpha} \rightarrow D \; \in \; Dov(D)$} be a collection of 
families of arrows (to be considered as basic covers on some 
basic objects D).

\vspace{1ex}
Define $Cov_{0} \;=\; Iso \;\cup\; Dov$, and:  

$Cov_{1} = \pi Cov_{0}$.

for an ordinal $\rho + 1$, 
$Cov_{\rho + 1} = Cov_{\rho} \circ Cov_{1}$.

for a limit ordinal $\rho$, 
$Cov_{\rho} = \bigcup _{\nu < \rho} Cov_{\nu}$.

\vspace{1ex}
Then 

1) $ \forall \, \nu < \rho \;\; Cov_{\nu} \subset Cov_{\rho}$.

2) $ \forall \, \rho \;\; Cov_{\rho}$ is a covering system.

3) $ \forall \, \rho \,,\; \nu  \;\;
                Cov_{\rho} \circ Cov_{\nu} \subset  Cov_{\rho + \nu}$ 
                (actually, equality holds).

\end{proposition}
\begin{proof}
1) is clear, 2) follows easily by induction. 3) follows for each 
$\rho$ by induction on $\nu \,$; on $\nu + 1$ by associativity, and on 
limit ordinals it is straightforward.    
\end{proof}

\begin{proposition} \label{generated}
With the notation in the previous proposition, 
 $Cov = \bigcup _{all \; \rho} Cov_{\rho}$ is the (pre) topology 
 generated by $Dov$.
\end{proposition}
\begin{proof}
It is clearly a covering system by \ref{pregenerated},2). 
It remains to see it is closed under composition. Consider 
$A_{\alpha, \, \beta} \rightarrow A_{\alpha} \rightarrow A$, with
\mbox{$A_{\alpha} \rightarrow A \; \in \; Cov_{\nu}(A)$}, and 
  $A_{\alpha,\, \beta} \rightarrow A_{\alpha} \; \in \; 
  Cov_{\rho_{\alpha}}(A_{\alpha})$. 
  Take an ordinal $\rho$ such that 
  \mbox{$\rho \geq \rho_{\alpha} \;\; \forall \, \alpha$} and use
  \ref{pregenerated},1) and \ref{pregenerated},3). 
 \end{proof}

\vspace{2ex}

\section{Proof of the theorems \ref{transitivity2} and \ref{Galois2}} 
\section*{in this section we assume the validity of \ref{assumption2}}

In \cite{G2}, Expose I, 10.2 and 10.3 Grothendieck defines the 
important concept of \emph{strict epimorphism}. Since then a whole variety of 
equivalent and/or related versions of this notion appeared in the literature under 
all sorts of names. To avoid confusion and fix the notation we recall now 
this original definition.

Let $X \mr{f} Y$ be an arrow in a category $\mathcal{C}$ and let $Ker_{f}$ 
be the full subcategory of the appropriate slice category whose objects are 
pairs of arrows
$\diagrama{ C \ar@<1ex>[r]^{x} \ar@<-1ex>[r]^{y} & X }$ such that 
$fx = fy$. Then: 

\begin{definition} \label{strictepi}
$f$ is an \emph{strict epimorphism} if for any other 
arrow  $X \mr{g} Z$ such that $Ker_{f} \subset Ker_{g}$, there exists a 
unique $Y \mr{h} Z$ such that $g = hf$. When $Ker_{f}$ has a terminal 
object the strict epimorphism is called \emph{effective}.   
\end{definition}

It immediately follows that strict epimorphisms are epimorphisms and 
that \emph{strict epi + mono = iso}.


\begin{proposition}  \label{poset}
The diagram of $F$, $\Gamma_{F}$, is a cofiltered category, and it is already 
a poset, $\Gamma_{F} \,=\, D_{F}$. 
\end{proposition}
\begin{proof}
By definition of prorepresentable functor, $\Gamma_{F}$ is a cofiltered 
category. Since all maps $X \mr{f} Y$ in $\mathcal{C}$ are epimorphisms,
for any object $Z \in \mathcal{C}$ the transition morphisms corresponding
to an arrow \mbox{$(x,X) \mr{f} (y,Y)$} in $\Gamma_{F}$,
$f^*:[Y,Z] \mr{} [X,Z]$ are all injective functions.
By construction of filtered colimits in $\mathcal{E}ns$ it follows
that the canonical maps of the colimit $[X,-] \mr{x} F$
are injective natural transformations (thus monomorphisms
in the category $\mathcal{E}ns^\mathcal{C}$). This implies that $\Gamma_{F}$ 
is a poset.
\end{proof}

\begin{proposition}  \label{refiso}
The functor $F$ is faithful and reflects isomorphisms.
\end{proposition}
\begin{proof}
Let $X \mr{f} Y$ be such that $F(f)$ is an isomorphism. We shall see 
that $f$ is a monomorphism (and thus, by \ref{assumption2} i), an 
isomorphism). Let $\diagrama{ Z \ar@<1ex>[r]^{t} \ar@<-1ex>[r]^{s} & X }$ 
be such that $fs = ft$. Clearly it follows that $F(s) = F(t)$. Take 
any $z \in FZ$ (use \ref{assumption2} ii) and let $x = F(s)(z) = 
F(t)(z)$. In this way $s$ and $t$ define arrows $(Z, z) \mr{} (X, x)$ 
in $\Gamma_{F}$. It follows from \ref{poset} that we must have $s = t$. 
Observe that within this argument we have also shown that $F$ is faithful.      
\end{proof}

\vspace{2ex}

Recall now the construction \ref{lAut(F)} of the locale $lAut(F)$.

\begin{definition}
Given an object  $[A]$ (determined by a 
finite subset $A \subset D_{\Delta F}$) on the site of 
definition $\mathcal{D}(D_{\Delta F})$, the \emph{content} of $[A]$
is the set of generators which are below $[A]$. That is, it is the set of 
objects $[(M, \; <\! m_{0}\,|\, m_{1} \!>)]$, 
$M \in \mathcal{C}$, $(m_{0},\; m_{1}) \in FM \times FM$ such that 
$[(M, \; <\! m_{0}\,|\, m_{1} \!>)] \, \leq \, [A]$. Notice that this means  
 \mbox{$[(M, \; <\!m_{0}\,|\, m_{1}\!>)] \:\leq\:[(X, \; <\!x_{0}\,|\, 
 x_{1}\!>)]$} for each 
 $(X,\; (x_{0},\; x_{1})) \in A$, 
 which in turn means that there is an arrow $M \mr{f} X$ 
 in $\mathcal{C}$ such that $f(m_{0}) = x_{0}\,,\; f(m_{1}) = x_{1}$ 
 (see \ref{freeinflattice} and \ref{deF}). 
\end{definition}
 
 \begin{proposition} \label{emptycontent}
 For each $X \in \mathcal{C}$, and each $x \neq y,\, z \in FX$,  
 \mbox{$[(X,\;<\!x\,|\,z\!>),\;(X,\;<\!y\,|\,z\!>)]$} and 
 $[(X,\;<\!z\,|\,x\!>),\;(X,\;<\!z\,|\,y\!>)]$ in 
 $\mathcal{D}(D_{\Delta F})$ have empty content.
 \end{proposition}
 \begin{proof}
 Let $[(M, \; <\!m_{0}\,|\, m_{1}\!>)] \:\leq\:
                                   [(X,\;<\!x\,|\,z\!>),\;(X,\;<\!y\,|\,z\!>)]$.
 Then, there are arrows $M \mr{f} X$ and $M \mr{g} X$ such that 
 \mbox{$f(m_{0}) = x$}, $f(m_{1}) = z\,,\;g(m_{0}) = y\,,\; g(m_{1}) = z$. 
 It follows from \ref{poset} that $f = g$, thus $x = y$, contrary 
 with the assumption. In the second case we do in the same way.                                    
 \end{proof} 

\begin{proposition} \label{auxiliar}
Let $W \mr{f} X$ and $W \mr{g} Y$ be any two arrows in $\mathcal{C}$, 
and $x_{0} \in FX$, $y_{0} \in FY$. If for $w \in FW$ the implication
\mbox{$F(f)(w) = x_{0}  \; \Rightarrow \; F(g)(w) = y_{0}$} holds, then there 
exists a unique $X \mr{h} Y$ such that $g = hf$ and $F(h)(x_{0}) = y_{0}$.
\end{proposition}


\begin{proof}
We prove first that under the assumption in the proposition, 
for arbitrary $v \in FW, \; w \in FW$, the 
following implication holds:
$$ 1) \;\; F(f)(v) = F(f)(w) \;\Rightarrow \; F(g)(v) = F(g)(w)$$

Take $M,\; m \in FM$ and $M \mr{s} W$, $M \mr{t} W$ such that 
$F(s)(m) = v$, $F(t)(m) = w$ (recall $\Gamma_{F}$ is cofiltered). It 
follows that \mbox{$F(fs)(m) = F(ft)(m)$,} thus by \ref{poset} $fs = ft$. 
By \ref{assumption2} i), iii)  take $m_{0} \in FM$ such that 
$F(fs)(m_{0}) = F(ft)(m_{0}) = x_{0}$. 
Let $v_{0} = F(s)(m_{0}), \; w_{0} = F(t)(m_{0})$. Clearly, 
$$[(M,\;<\!m,\;m_{0}\!>)] \rightarrow 
[(W,\;<\!v,\;v_{0}\!>),\; (W,\;<\!w,\;w_{0}\!>)]$$ and 
$F(f)(v_{0}) = F(f)(w_{0}) = x_{0}$. We have also 
$$[(W,\;<\!v,\;v_{0}\!>)] \rightarrow 
[(Y,\;<\!F(g)(v),\;F(g)(v_{0})\!>)]$$  
$$[(W,\;<\!w,\;w_{0}\!>)] \rightarrow 
[(Y,\;<\!F(g)(w),\;F(g)(w_{0})\!>)]$$
By assumption $F(g)(v_{0}) = F(g)(w_{0}) = y_{0}$. Thus, we have 
$$[(M,\;<\!m,\;m_{0}\!>)] \rightarrow 
[(Y,\;<\!F(g)(v),\; y_{0})\!>),\; (Y,\;<\!F(g)(w),\; y_{0})\!>)]$$ 
It follows then from \ref{emptycontent} that we must have  
$F(g)(v) = F(g)(w)$. This finishes the proof of 1).

It follows from 1) that $Ker_{F(f)} \subset Ker_{F(g)}$. Since F 
is faithful (\ref{refiso}) this implies $Ker_{f} \subset Ker_{g}$. 
The proof finishes by definition of strict epimorphism 
(\ref{strictepi}).      
\end{proof}

\vspace{2ex}

Our first important result is the following:

\begin{theorem} \label{importante1}
Let $[A_{\alpha}] \rightarrow [A]$ in  $\mathcal{D}(D_{\Delta F})$ be any cover 
in the site of definition 
of $lAut(F)$ (see \ref{lAut(F)}). Then, if $[A]$ has nonempty 
content, there exists an index $\alpha$ such that $[A_{\alpha}]$ has non 
empty content. 
\end{theorem}
\begin{proof}
By induction on the generation of covers (see  
\ref{pregenerated} and \ref{generated}).

\vspace{1ex}
$1$) Let [D] be the object in the two basic empty covers: \newline
\mbox{$[D] = [(X,\;<\!x\,|\,z\!>),\;(X,\;<\!y\,|\,z\!>)]$ or 
$[(X,\;<\!z\,|\,x\!>),\;(X,\;<\!z\,|\,y\!>)]$}, 
By \ref{emptycontent} a pullback of 
the form
$$
\diagrama
     {
       \emptyset  \ar[r] \ar[d]  &  [A] \ar[d]  \\
       \emptyset  \ar[r]  &  [D] 
     }
$$
can not be, since it implies that $[D]$ would have non 
empty content.


Consider now $[D] = 1$, $Z \in \mathcal{C}$, $z_{1} \in FZ$, the basic cover  
\mbox{$[(Z,\;<\!z\,|\,z_{1}\!>)] \;\rightarrow \; 1, \; z \in FZ$}, and the 
pullback (see \ref{freeinflattice2}):
$$
\diagrama
     {
       [(Z,\,<\!z\,|\,z_{1}\!>),\; A]  \ar[r] \ar[d]  &  [A] \ar[d]  \\
       [(Z,\,<\!z\,|\,z_{1}\!>)]  \ar[r]  &  1 
     }
$$
Let  
$[(M, \; <\! m_{0}\,|\, m_{1} \!>)] \, \leq \, [A]$, and take 
$(N,\; n_{1}) \rightarrow (Z,\; z_{1})$,
$(N,\; n_{1}) \rightarrow (M,\; m_{1})$ in $\Gamma_{F}$. Since the 
function $FN \rightarrow FM$ is surjective (\ref{assumption2}), we can take 
$n_{0} \in FN$ such that $n_{0} \mapsto m_{0}$, and let $z_{0}$ be the 
image of $n_{0}$ in $FZ$, $n_{0} \mapsto z_{0}$. We have then
\mbox{$[(N, \; <\! n_{0}\,|\, n_{1} \!>)] \, \leq \,  
[(Z,\,<\!z_{0}\,|\,z_{1}\!>)\,,\; (M, \; <\! m_{0}\,|\, m_{1} \!>)]$}.  
This shows that $[(Z,\,<\!z_{0}\,|\,z_{1}\!>),\; A]$ 
(corresponding to the index $z_{0}$ in the cover) has non empty content.

The same argument applies to the remaining basic covers 
\mbox{$[(Z,\;<\!z_{0}\,|\,z\!>)] \;\rightarrow \; 1, \; z \in FZ$}.

\vspace{1ex}
$\rho + 1$) Consider now the cover 
$[A_{\alpha, \, \beta}] \rightarrow [A_{\alpha}] \rightarrow [A]$, with 
\mbox{$[A_{\alpha, \, \beta}] \rightarrow [A_{\alpha}]\; \in \; 
Cov_{\rho}$}  and
$[A_{\alpha}] \rightarrow [A] \; \in \; Cov_{1}$. Take $\alpha$ such that
$[A_{\alpha}]$ has non empty content, and for this $\alpha$ take 
$\beta$ such that $[A_{\alpha, \, \beta}]$ has non empty content.

\vspace{1ex}
limit ordinal $\rho$) In this case the proof is even more immediate.
\end{proof}

\begin{corollary} \label{corolarioimportante1}
If $[A] \in \mathcal{D}(D_{\Delta F})$ has non empty content, then the 
empty family does not cover $[A]$. In particular, for any 
\mbox{$X \in \mathcal{C}$}, $(x_{0},\; x_{1}) \in FX \times FX$, the empty 
family does not cover \mbox{$[(X, \; <\! x_{0}\,|\, x_{1} \!>)]$}.
\end{corollary}
 \begin{proof}
Clear, since for the empty cover can not exist any index.
\end{proof}  

The fact that the empty family does not cover 
\mbox{$[(X, \; <\! x_{0}\,|\, x_{1} \!>)]$} means that this object 
stays different from $0$ in the sheaf poset.

This proves theorem \ref{transitivity2}.
 
\begin{corollary} [Theorem \ref{transitivity2}] \label{transitivity3}
For each $X \in \mathcal{C}$, the action of $lAut(F)$ on the set $FX$ 
(defined in \ref{actionF}) is transitive. Explicitly, 
\mbox{$\forall\: (x_{0},\; x_{1}) \in FX \times FX$,}  
$\#[(X,\,<\!x_{0}\,|\,x_{1}\!>)] \:\neq \: 0$ (see \ref{action}).
\end{corollary}

\vspace{2ex}

Our second important result is the following:

\begin{theorem} \label{importante2}
Given two objects 
$[(X,\,<\!x_{0}\,|\,x_{1}\!>)]$, $[(Y,\,<\!y_{0}\,|\,y_{1}\!>)]$ 
in $D_{\Delta F}$,
and a cover\mbox{$[A_{\alpha}] \rightarrow [(X,\,<\!x_{0}\,|\,x_{1}\!>)]$} 
in $\mathcal{D}(D_{\Delta F})$ in the site of definition of $lAut(F)$ (see \ref{lAut(F)}), 
 the following implication holds:
$$
 \forall\, \alpha \; [A_{\alpha}] \rightarrow  [(Y,\,<\!y_{0}\,|\,y_{1}\!>)]
\;\; \Longrightarrow \;\; [(X,\,<\!x_{0}\,|\,x_{1}\!>)] \rightarrow 
  [(Y,\,<\!y_{0}\,|\,y_{1}\!>)]
$$          
\end{theorem}
\begin{proof}
Consider first that by the corollary above a cover 
of \mbox{$[(X,\,<\!x_{0}\,|\,x_{1}\!>)]$} can not be empty.
Let now \mbox{$[(Z,\;<\!z\,|\,z_{1}\!>)] \;\rightarrow \; 1, \; z \in FZ$} 
be one of the other basic covers, and consider the  $Cov_{1}$ cover 
determined by following pull-back:
$$
\diagrama
     {
       [(Z,\,<\!z\,|\,z_{1}\!>),\; (X,\,<\!x_{0}\,|\,x_{1}\!>)]
                \ar[r] \ar[d]  &  [(X,\,<\!x_{0}\,|\,x_{1}\!>)] \ar[d]  \\
       [(Z,\,<\!z\,|\,z_{1}\!>)]  \ar[r]  &  1 
     }
$$


Take 
$(M,\; m_{1}) \rightarrow (Z,\; z_{1})$,
$(M,\; m_{1}) \rightarrow (X,\; x_{1})$ in $\Gamma_{F}$.
Consider all the  $m \in FM$ such that $m \mapsto x_{0}$, and let $z$ be the 
images of these $m$ in $FZ$, $m \mapsto z$. This defines, for each 
such $m$   
$$
[(M,\,<\!m\,|\,m_{1}\!>)] \rightarrow 
  [(Z,\,<\!z\,|\,z_{1}\!>),\; (X,\,<\!x_{0}\,|\,x_{1}\!>)].
$$
  
We start now the induction on the generation of covers 
(see \ref{pregenerated} and \ref{generated}). 
We deal simultaneously with the case $\rho = 1$ and the case $\rho + 1$.
Consider the $Cov_{\rho + 1}$ cover determined by a $Cov_{1}$ cover as above,
 and for each $z \in FZ$, a $Cov_{\rho}$ cover 
$$[A_{z,\;\alpha}] \rightarrow 
 [(Z,\,<\!z\,|\,z_{1}\!>),\; (X,\,<\!x_{0}\,|\,x_{1}\!>)]$$ 
(the case $\rho = 1$ is included considering all these covers to be the 
identity). 

For each $m$ (and $z$, $m \mapsto z$) as above, consider the following diagram, 
defined as a pullback in $\mathcal{D}(D_{\Delta F})$
$$
\diagrama
     {
      [B_{z,\;\alpha}] \ar[r] \ar[d] & [(M,\,<\!m\,|\,m_{1}\!>)] \ar[d] \\
      [A_{z,\;\alpha}] \ar[r] & 
      [(Z,\,<\!z\,|\,z_{1}\!>),\; (X,\,<\!x_{0}\,|\,x_{1}\!>)]   
     }
$$
By assumption, for each $\alpha$, there is
$[A_{z,\;\alpha}] \rightarrow [(Y,\,<\!y_{0}\,|\,y_{1}\!>)]$. Composing 
we have $[B_{z,\;\alpha}] \rightarrow [(Y,\,<\!y_{0}\,|\,y_{1}\!>)]$. 
Finally, by
 \ref{pregenerated}, 2) and the inductive hypothesis we have 
$[(M,\,<\!m\,|\,m_{1}\!>)] \rightarrow [(Y,\,<\!y_{0}\,|\,y_{1}\!>)].$
 Since all these arrows (one for each $m$) send $m_{1} \mapsto y_{1}$, 
 by \ref{poset} they all correspond to a 
same single arrow $M \rightarrow Y$ in $\mathcal{C}$.  

In conclusion, we have two arrows $M \rightarrow X$, $M \rightarrow Y$ such 
that for $m \in FM$, if \mbox{$m \mapsto x_{0}$}, then  $m \mapsto y_{0}$. 
It follows by \ref{auxiliar} that there exist 
$X \rightarrow Y$ such that $x_{0} \mapsto y_{0}$. Since the composite 
$M \rightarrow X \rightarrow Y$ is the arrow $M \rightarrow Y$, 
it is also the case that $x_{1} \mapsto y_{1}$. Thus we have
\mbox{$[(X,\,<\!x_{0}\,|\,x_{1}\!>)] \rightarrow 
   [(Y,\,<\!y_{0}\,|\,y_{1}\!>)]$.}
   
The same argument applies to the other remaining basic covers 
\mbox{$[(Z,\;<\!z_{0}\,|\,z\!>)] \;\rightarrow \; 1, \; z \in FZ$}.
     
The case of a limit ordinal is evident. This finishes the proof    
\end{proof}

\vspace{2ex}

Clearly the topology on $\mathcal{D}(D_{\Delta F})$ that defines 
$lAut(F)$ is not subcanonical, and so the morphism of inf-posets 
$\mathcal{D}(D_{\Delta F}) \mr{\#} lAut(F)$ (where $\#$ indicates 
the associated sheaf) is far from being full. However, for the full subposet 
$D_{\Delta F} \hookrightarrow \mathcal{D}(D_{\Delta F})$ the theorem 
above gives:

\begin{corollary} \label{full}
The morphism of posets $D_{\Delta F} \mr{\#} lAut(F)$ is full. 
Explicitly, if 
$\#[(X,\,<\!x_{0}\,|\,x_{1}\!>)] \rightarrow \#[(Y,\,<\!y_{0}\,|\,y_{1}\!>)]$ 
in $lAut(F)$, then there exists a unique $X \rightarrow Y$ in 
$\mathcal{C}$ such that $x_{0} \mapsto y_{0}$ and $x_{1} \mapsto y_{1}$. 
\end{corollary}


\begin{proof}
Consider the following chain of equivalences (or bijections) 
justified, in turn, by definition of $\#$, (Yoneda and) construction 
of $\#$, and theorem \ref{importante2} respectively:

\vspace{1ex} 

$\begin{array}{c}
\#[(X,\,<\!x_{0}\,|\,x_{1}\!>)] \rightarrow 
\#[(Y,\,<\!y_{0}\,|\,y_{1}\!>)]  \\[1ex] \hline
[(X,\,<\!x_{0}\,|\,x_{1}\!>)] \rightarrow 
\#[(Y,\,<\!y_{0}\,|\,y_{1}\!>)] \\[1ex] \hline 
\exists \;cover \; [A_{\alpha}] \rightarrow [(X,\,<\!x_{0}\,|\,x_{1}\!>)] \;| 
\;\;  \forall\, \alpha \; [A_{\alpha}] \rightarrow  
[(Y,\,<\!y_{0}\,|\,y_{1}\!>)] \; \\[1ex] \hline 
[(X,\,<\!x_{0}\,|\,x_{1}\!>)] \rightarrow 
[(Y,\,<\!y_{0}\,|\,y_{1}\!>)]    \\
\end{array}$

\end{proof}

This proves theorem \ref{Galois2}.

\begin{corollary} [Theorem \ref{Galois2}, Lifting Lemma]
Given any objects 
$X \in \mathcal{C},\, Y \in \mathcal{C},\, and \, x \in FX,\, y \in FY$, 
if \mbox{$lFix(x) \leq lFix(y)$} in $lAut(F)$, then there exist a unique 
arrow \mbox{$X \mr{f} Y$} in $\mathcal{C}$ such that $F(f)(x) = y$. 
\end{corollary}
\begin{proof}
Notice that the Galois group 
$lFix(x)$ for the action of $lAut(F)$ on $FX$ is given by 
$lFix(x) \;=\; \#[(X,\,<\!x\,|\,x\!>)]$ (see \ref{lFix(x)}). Thus, clearly, 
this statement is the particular case of \ref{full}, when $x_{0} = x_{1}$ 
and \mbox{$y_{0} = y_{1}$}.
\end{proof} 

\vspace{2ex}

\section{preliminaries on the classifying topos of a localic 
group}

Given a set $X$, by the construction in proposition \ref{lAut(X)}, 
the following equations hold in the local $lAut(X)$ (recall that we abuse the 
notation and omit to indicate the associate sheaf morphism):
$$[<\!z\,|\,x\!>,\;<\!z\,|\,y\!>] \;=\; 0\,,\;\;\; 
  [<\!x\,|\,z\!>,\;<\!y\,|\,z\!>] \;=\; 0 \;\;\; 
(each \; x \neq y,\; z)$$  
$$\bigvee\nolimits_{x} \; [<\!x\,|\,z\!>] \;=\; 1\,, \;\;\;
  \bigvee\nolimits_{x} \; [<\!z\,|\,x\!>] \;=\; 1 \;\;\;
(each \, z)$$

Recall also that a \emph{morphism of localic groups} $H \mr{\varphi} G$ is a 
continuous map such that 
$m^{*}\varphi^{*} = (\varphi^{*}\otimes \varphi^{*})\,m^{*}$
(where $m$ denotes the multiplication in the two structures).  
 
\begin{definition} \label{action}

Given a localic group $G$ and a set $X$, an \emph{action} of $G$ on 
$X$ is a continuous morphism of localic groups $G \mr{\mu} lAut(X)$. It 
is completely determined by the value of its inverse image on the 
generators, $X \times X  \mr{\mu^{*}} G$. We say that the action is 
\emph{transitive} 
when for all \mbox{$x \in X,\: y \in X$, $\mu^{*}[<x\,|\,y>] \neq 0$}.  
\end{definition}

\begin{definition} \label{lFix(x)}
Given a localic group $G$ acting on a set $X$, and element $x\in X$, 
the open subgroup of $G$, described informally
 as $\{g\in G \;|\; gx = x \}$, is defined to be the object 
 \mbox{$lFix(x) = \mu^{*}[<x\,|\,x>]$} in the locale $G$.
\end{definition}


Given a localic group $G$, a \emph{$G$-set} is a set furnished with 
an action of $G$. Given two $G$-sets $X$, $Y$,
 $X \times X  \mr{\mu^{*}} G$, \mbox{$Y \times Y  \mr{\mu^{*}} G\,$,} a morphism 
 of $G$-sets is a function $X \mr{f} Y$ such that 
\mbox{$\mu^{*}[<\!x \:|\:y\!>] \;\leq\; \mu^{*}[<\!f(x) \:|\:f(y)\!>]$.} 
This defines a category $\mathcal{B}G$ furnished with an  
underline set functor $\mathcal{B}G \mr{|\;\;|} \mathcal{E}ns$ 
into the category of sets. We shall denote 
$t\mathcal{B}G$ the full subcategory of non empty transitive $G$-sets.

It is easy to check the following (consider \ref{lAut(F)} and 
\ref{actionF}):

\begin{proposition} \label{GtoF}
Let $F$ be the underline set functor $\mathcal{B}G \mr{|\;\;|} \mathcal{E}ns$.
 The map given by $[(X,\,<\!x_{0}\,|\,x_{1}\!>)] \longmapsto
 \mu^{*}[<\!x_{0} \:|\:x_{1}\!>]$ determines (the inverse image of) a 
 morphism of localic groups $G \rightarrow lAut(F)$.  
\end{proposition} 

\begin{proposition} \label{suryectivos}
All morphisms $X \mr{f} Y$ between non empty transitive $G$-sets are surjective 
functions of the underline sets.
\end{proposition}
\begin{proof}
We shall see that 
$ \forall \, y \in Y \;\; \exists \, x \in X \;\;|\;\; f(x) = y$.

Take any $x_{0} \in X$ and let $y_{0} = f(x_{0})$. Then:
$$ 
1 \;=\; \bigvee\nolimits_{x}\, \mu^{*}[<\!x_{0} \:|\: x\!>] \;\leq\;
\bigvee\nolimits_{x}\, \mu^{*}[<\!y_{0} \:|\: f(x)\!>]$$
Taking the infimum against $\mu^{*}[<\!y_{0} \:|\: y\!>]$,
$$ 0 \;\neq\; \mu^{*}[<\!y_{0} \:|\: y\!>] \;\leq\; 
\bigvee\nolimits_{x}\, \mu^{*}[<\!y_{0} \:|\: y\!>,\; <\!y_{0} \:|\: f(x)\!>]$$
The terms in the supremum are equal to $0$ except if $y = f(x)$. This 
finishes the proof.
\end{proof}   

\begin{proposition} \label{descomposicion}
Given a localic group acting on set, $G \mr{\mu} lAut(X)$,
$X \times X  \mr{\mu^{*}} G$, the relation 
\mbox{$x \sim y \; \Leftrightarrow \; \mu^{*}([<\!x \:|\:y\!>]) \neq 0$} is an 
equivalence 
relation on $X$, and G acts transitively on each equivalence class.
\end{proposition}
\begin{proof}
Assume $\mu^{*}[<\!x \:|\:z\!>] \neq 0$ and $\mu^{*}[<\!z \:|\:y\!>] \neq 0$.
The multiplication $m$ of $lAut(X)$ is given by 
$$m^{*}([<\!x \:|\:y\!>]) \;=\;
 \bigvee\nolimits_{z} \;[<\!x \:|\:z\!>]\otimes[<\!z \:|\:y\!>]\;\;\;$$
Since $\mu^{*}$ is a morphism of groups as well as of locales, 
$$m^{*}\mu^{*}[<\!x \:|\:y\!>] \;=\;
 (\mu^{*}\otimes \mu^{*})\,m^{*}[<\!x \:|\:y\!>] \;=\; 
\bigvee\nolimits_{z} \;\mu^{*}[<\!x \:|\:z\!>]\otimes\mu^{*}[<\!z \:|\:y\!>]$$ 
It follows  that $m^{*}\mu^{*}([<\!x \:|\:y\!>]) \neq 0$. Thus 
$\mu^{*}([<\!x \:|\:z\!>]) \neq 0$. The second assertion is obvious.    
\end{proof}
    

Given an element $x_{0} \in G$, the \emph{connected component} of $x_{0}$ is 
the transitive $G$-set with underline set 
$\{x \in X\;|\;\mu^{*}([<\!x \:|\:x_{0}\!>]) \neq 0\}$. 

The coproduct of $G$-sets is just the disjoint union furnished with 
the obvious action ($\mu^{*}[<\!x \:|\:y\!>] = 0$ if $x$ and $y$ are 
in different components). In this way, every $G$-set is the coproduct 
in $\mathcal{B}G$ of transitive $G$-sets. With this it is clear that 
it follows from \ref{descomposicion} that \emph{a $G$-set is a connected  object
in $\mathcal{B}G$ if and only if the action is transitive}.


\begin{proposition} \label{cofiltrante}
The diagram of the underline set functor 
\mbox{$t\mathcal{B}G \mr{|\;\;|} \mathcal{E}ns$} (from the category of 
non empty transitive $G$-sets) is a cofiltered poset. That 
is:

1) Given morphisms of transitive $G$ sets,
$\diagrama{ X \ar@<1ex>[r]^{t} \ar@<-1ex>[r]^{s} & Y }$ and 
$x_{0} \in X$ such that $s(x_{0}) = t(x_{0})$, then s = t.

2) Given two transitive $G$-sets $X$, $Y$, and elements $x \in X$, 
$y \in Y$, there exists a transitive $G$-set $M$, an element $m \in M$,
and morphisms of $G$-sets $M \mr{s} X$, $M \mr{t} Y$, such 
that $s(m) = x$, $t(m) = y$. 
\end{proposition}
\begin{proof}

1) Let $y_{0} = s(x_{0}) = t(x_{0})$, and let $x$ be any element in $X$. 
Since $0 \neq \mu^{*}([<\!x_{0} \:|\:x\!>]$ it follows that 
$0 \neq \mu^{*}[<\!y_{0} \:|\: s(x)\!>,\; <\!y_{0} \:|\: t(x)\!>]$. 
Thus it must be $s(x) = t(x)$
 
2) Take the connected component of $(x,\:y)$ in the product 
$X \times Y$ and the two projections (the action in the product is given by
\mbox{$\mu^{*}[<\!(x,\; y)\:|\:(x',\;y')\!>] = 
\mu^{*}[<\!x \:|\:x'\!>] \wedge  \mu^{*}[<\!y \:|\:y'\!>]$)}.  
\end{proof}

\begin{proposition} \label{transstrictepi}
All morphisms $X \mr{f} Y$ between non empty transitive $G$-sets are 
strict epimorphisms in $t\mathcal{B}G$.
\end{proposition}
\begin{proof}
Let $X \mr{g} Z$ be such that $Ker_{f} \subset Ker_{g}$ (see 
\ref{strictepi}).

 We shall see 
first that $Ker_{|\,f\,|} \subset Ker_{|\,g\,|}$ taken in 
$\mathcal{E}ns$. Let $x \in X$, $y \in X$ be such that $f(x) = f(y)$. Take 
$M$, $m \in M$, $s$ and $t$ as in \ref{cofiltrante},2). Then $fs(m) = ft(m)$, 
and thus $s = t$. By assumption it follows $gs = gt$, thus $g(x) = g(y)$. 

From this, since $f$ is surjective (\ref{suryectivos}), it follows 
there exists a function $Y \mr{h} Z$ such that $hf = g$. It remains to 
see that $h$ is a morphism of $G$-sets. We do this now.

Let $y_{0},\; y_{1}$ be any two points in Y. Take $x_{0} \in X$, 
$f(x_{0}) = y_{0}$. We have
$\mu^{*}[<\!y_{0} \:|\: y_{1}\!>] \wedge \mu^{*}[<\!x_{0} \:|\: x\!>] 
\; \leq \; \mu^{*}[<\!y_{0} \:|\: y_{1}\!>,\;<\!y_{0} \:|\: 
f(x)\!>]$, which equals $0$ unless $f(x) = y_{1}$. With this:
$$\mu^{*}[<\!y_{0} \:|\: y_{1}\!>] \;=\; \bigvee\nolimits_{x}\;
\mu^{*}[<\!y_{0} \:|\: y_{1}\!>] \wedge \mu^{*}[<\!x_{0} \:|\: x\!>] 
\;=$$
$$\bigvee\nolimits_{f(x) \,=\, y_{0}}\;
\mu^{*}[<\!y_{0} \:|\: y_{1}\!>] \wedge \mu^{*}[<\!x_{0} \:|\: x\!>]
 \;\leq\;
 \bigvee\nolimits_{f(x) \,=\, y_{0}}\;\mu^{*}[<\!x_{0} \:|\: x\!>] \;\leq$$
$$ \bigvee\nolimits_{f(x) \,=\, y_{0}}\;\mu^{*}[<\!g(x_{0}) \:|\: g(x)\!>]
 \;=\; \mu^{*}[<\!h(y_{0}) \:|\: h(y_{1})\!>]$$        
\end{proof}
    

Clearly, by definition, given a morphism of $G$-sets $X \mr{h} Y$, if 
\mbox{$x_{0} \in X$,} and $y_{0} = h(x_{0})$, we have 
$lFix(x_{0}) \leq lFix(y_{0})$ in $G$. The reverse implication also 
holds, which means that transitive $G$-set are in a sense quotients of $G$. 


\begin{proposition} \label{cociente}
Let $X$ be any transitive $G$-set, $x_{0} \in X$. Given any $G$-set $Y$,
 $y_{0} \in Y$, such that $lFix(x_{0}) \leq  lFix(y_{0})$ in G, there 
 exists a unique morphism $X \mr{h} Y$ such that $h(x_{0}) = y_{0}$.
\end{proposition}
\begin{proof}
Take $M$, $m_{0} \in M$, $M \mr{f} X$ and $M \mr{g} Y$, 
$f(m_{0}) = x_{0}$, $g(m_{0}) = y_{0}$ (\ref{cofiltrante},2).
To prove the statement it is enough to show that
 $Ker_{f} \subset Ker_{g}$ (\ref{transstrictepi}). First we prove 
the following implication:
$$\;\;\; 1)  \;\;\;\forall\,m \in M,\;\;  
         f(m) = x_{0} \; \Rightarrow \; g(m) = y_{0}$$
Assume $f(m) = x_{0}$. Then,
$\mu^{*}[<\!m_{0} \:|\: m\!>] \leq \mu^{*}[<\!y_{0} \:|\: g(m)\!>]$, 
and,  
$ \mu^{*}[<\!m_{0} \:|\: m\!>] \leq \mu^{*}[<\!x_{0} \:|\: x_{0}\!>] 
\leq \mu^{*}[<\!y_{0} \:|\: y_{0}\!>]$.
Thus \mbox{$0 \leq \mu^{*}[<\!y_{0} \:|\: y_{0}\!>,\;
<\!y_{0} \:|\: g(m)\!>]$.} It follows $g(m) = y_{0}$.

With this, now we prove
 $Ker_{f} \subset Ker_{g}$.
 Let $\diagrama{ Z \ar@<1ex>[r]^{s} \ar@<-1ex>[r]^{t} & Y }$ be in 
 $Ker_{f}$, that is, $fs = ft$. Take $z_{0}$ such that 
 $fs(z_{0}) = ft(z_{0}) = x_{0}$ (\ref{suryectivos}). Then, for any 
 $z \in Z$,
$$0 \;\leq\; \mu^{*}[<\!z \:|\: z_{0}\!>] \;\leq\; 
\mu^{*}[<\!gs(z) \:|\: gs(z_{0})\!>] \wedge 
\mu^{*}[<\!gt(z) \:|\: gt(z_{0})\!>]$$ 
$$\;\;\;\;\;\;\;\;\;\;\;\leq
\mu^{*}[<\!gs(z) \:|\: y_{0}\!>,\; <\!gt(z) \:|\: y_{0}\!>]\,,$$
 the last inequality justified by 1). It follows that $gs(z) = gt(z)$, thus 
$gs = gt$. This finishes the proof.
\end{proof}

\begin{corollary} \label{small}
If the localic group $G$ is small (meaning it has only a set of 
objects), then the category $t\mathcal{B}G$ is also small.
\end{corollary}
\begin{proof}
Let $F$ be the underline set functor
 $t\mathcal{B}G \rightarrow \mathcal{E}ns$. Then the map
 $\Gamma_{F} \rightarrow G$ given by
 $(X, x_{0}) \longmapsto \mu^{*}[<\!x_{0} \:|\:x_{0}\!>]$ creates 
 (thus also reflects) isomorphisms (compare with \ref{GtoF}).  
\end{proof}

\vspace{2ex}
 
\section{Characterization of the classifying topos of a localic 
group}

In this section we characterize the category $\mathcal{B}G$ of 
$G$-sets in terms of the theory of topoi. That is, we prove Theorem B 
in the introduction. We shall see how this 
characterization follows in an straightforward manner from theorems   
\ref{Galois2} and \ref{transitivity2}.

First recall that a \emph{connected atomic topos} is a connected,
locally connected and boolean topos. The reference for atomic topoi and atomic 
sites is \cite{BD}. For connected and locally connected topoi see 
\cite{G2}, Expose IV, 2.7.5, 4.3.5, 7.6 and 8.7.

We have:

\begin{proposition} \label{atomicsite}
Let $\mathcal{C}$ be a category and $F:\mathcal{C} \mr{} 
\mathcal{E}ns$ be a functor as in \ref{assumption2}, i) and iii). Then the 
canonical (in this case atomic) topology defines an atomic site with 
a point. If $\mathcal{C}$ is small, the topos of sheaves 
$\mathcal{C^{\sim}}$ is an atomic topos with a point (see \cite{BD}).
This topos is connected if and only if condition ii) holds.
Any connected atomic topos with a point can be presented in this way. 
\end{proposition}

\begin{proposition} \label{gsetatomic}
The category $t\mathcal{B}G$ of transitive $G$-sets satisfies 
\ref{assumption2}. If $G$ is small, $\mathcal{B}G$ is an atomic topos 
with a point, with inverse image given by the underline set ($\mathcal{B}G$ 
is the topos of sheaves for the canonical topology on $t\mathcal{B}G$).  
\end{proposition}
\begin{proof}
The first assertion is given by \ref{transstrictepi} and 
\ref{suryectivos}. The second follows from \ref{descomposicion} and 
\ref{small}.
\end{proof}

Consider now any category  $\mathcal{C}$ and any set valued functor
$\mathcal{C} \mr{F} \mathcal{E}ns$. Then, proposition \ref{actionF} 
shows that $F$ lifts into a functor, that we denote $\mu F$, 
$\mathcal{C} \mr{\mu F} \mathcal{B}G$, for $G = lAut(F)$. We have: 

\begin{theorem} \label{joyaltierney}
Let $\mathcal{C}$ and $\mathcal{C} \mr{F} \mathcal{E}ns$ be as in 
\ref{assumption2}. Then the 
functor $\mu F$ lands into $t\mathcal{B}G$, 
$\mathcal{C} \mr{\mu F} t\mathcal{B}G$. If $\mathcal{C}$ is small, 
then $G$ is small, and $\mu F$ induces an equivalence 
of categories $\mathcal{C}^{\sim} \mr{\simeq} \mathcal{B}G$ between 
the topoi of sheaves for the canonical topology on $\mathcal{C}$ and 
the classifying topos $\mathcal{B}G$. 
\end{theorem}
\begin{proof}
Theorem \ref{transitivity2} just says that $\mu F$ lands into
$t\mathcal{B}G$. We shall prove:

1) The functor $\mu F$ (which is faithful since $F$ is, see \ref{refiso})
is also full.

2) Given any transitive $G$-set $S$, there exists $X \in \mathcal{C}$ 
and an strict epimorphism $\mu FX \rightarrow S$ in $t\mathcal{B}G$.

proof of 1). This is just the meaning of Theorem \ref{Galois2}. Given 
a morphism of $G$-sets $\mu FX \mr{f} \mu FY$, choose any 
$x_{0} \in FX$, and let $y_{0} = f(x_{0})$. By definition
$lFix(x_{0}) \leq lFix(y_{0})$ in $G$.

proof of 2). Choose any $s_{0} \in S$. Consider $lFix(s_{0}) \in G$. 
Clearly $e \in lFix(s_{0})$. Then, by 
\ref{identidad}, and the construction of $lAut(F)$ (\ref{lAut(F)}), since 
$\Gamma_{F}$ is cofiltered, it follows that there is $X \in \mathcal{C}$ and 
 $x_{0} \in FX$, such that 
 $lFix(x_{0}) \;=\; \#[(X,\,<\!x_{0}\,|\,x_{0}\!>)] \;\leq\; lFix(s_{0})$. The 
 proof finishes then by \ref{transstrictepi} and \ref{cociente}.
 
 The theorem follows from 1) and 2)  by the comparison lemma (\cite{G2}, Expose 
III, 4.).
 
\end{proof}

We can now easily establish Theorem B.

\begin{theorem} 
A topos $\mathcal{E}$
with a point $\mathcal{E}ns \mr{p} \mathcal{E}$, $p^{*} = F$, is 
connected atomic if and 
only if it is the classifying topos $\mathcal{B}G$ of a localic group $G$, 
and this group can be taken to be $lAut(F) = lAut(p)^{op}$.
\end{theorem}
\begin{proof}
From \ref{atomicsite} it follows that the easy direction on this 
equivalence is given by \ref{gsetatomic}, and the hard direction  
by \ref{joyaltierney}. 
\end{proof}

\end{document}